\documentstyle[11pt]{article}
\newtheorem {th}{Theorem}

\newtheorem {lem}[th]{Lemma}

\def\Cox{\hfill \Box}

\def\diseq{\, {\stackrel {{\cal D}} {=}}}

\def\sf{\sigma\mbox{-field}}

\def\E{{\bf{E}}}
\def\P{{\bf{P}}}
\def\N{\hbox{I\kern-.2em\hbox{N}}}
\def\R{\hbox{I\kern-.2em\hbox{R}}}

\def\Z{{\bf{Z}}}
\def\F{{\cal{F}}}
\def\G{{\cal{G}}}

\def\cE{{\cal{E}}}
\def\AP{A^p}
\def\DAP{\Delta A^p}

\def\AO{{}^o \kern-.3em A}
\def\DAO{\Delta \AO}

\def\|{\, | \, }
\def\v0{{\bf 0}}
\def\one{{\bf 1}}
\def\0{\hat{0}}
\def\1{\hat{1}}

\begin{document}
\begin{center}
{\large \bf MAXIMUM VARIATION OF TOTAL RISK}
\end{center}
\vspace{5ex}
\begin{flushright}
Robin Pemantle \footnote{Research supported in part by 
National Science Foundation grant \# DMS 9300191, by a Sloan Foundation
Fellowship, and by a Presidential Faculty Fellowship}$^,$\footnote{Department
of Mathematics, University of Wisconsin-Madison, Van Vleck Hall, 480 Lincoln
Drive, Madison, WI 53706}
  ~\\
\end{flushright}
\vspace{3ex}

\noindent{\large \bf Abstract:}  Let $Z > 0$ be a random time.  The total
risk of discovering $Z$ in the next time interval $(t , t+dt)$ is never
more variable than an exponential of mean one, which is achieved when
the information up to time $t$ is $\sigma (Z \wedge t)$.  \\[3ex]

\section{Results}

{\em Scenario 1:}  You have a life insurance policy for one million dollars.
The mortality tables for the entire population tell you that a lifespan of
$n$ years has probability $q_n$.  Your premium for year $n$ is 
$\$ 1,000,000 \cdot h_n$, where $h_n = q_n / \sum_{k=n}^\infty q_k$.
In the absence of further information this is fair: you may choose each
year whether to renew your policy, and your expected gain is always zero.
If further information becomes available each year, the fair premium
becomes $\$ 1,000,000 \cdot Q_n$, where $Q_n$ is the conditional 
probability of dying in year $n$ given all the information up to that point.
How does the extra information affect the distribution of the lifetime
total you pay for your policy?  

{\em Scenario 2:}  Random variables $\{ X(e) \}$
are assigned to the edges of a graph.  These values determine a random
subset $S$ of edges, called pivotal bonds.  You know that $|S| \leq K$
with high probability.  You order the edges $e(1) , e(2) , \ldots$ and
look at the values $X(e_j)$ one at a time.  You are interested (for reasons
explained below) in the distribution of the random variable 
$$W = \sum_{j=1}^\infty \P (e(j) \in A \| X(e(i)) : i < j) .$$
What bound can you get on $\P (W > \lambda K)$?  

The purpose of this note is to prove an inequality that answers the
questions in the two scenarios.  The relevant notion of variability turns
out to be the following one.  Define a partial order $\preceq$ among 
random variables by 
$$ Y \preceq X \; \mbox{ if and only if for every convex } \phi, 
   \E \phi (Y) \leq \E \phi (X) ,$$
where both expectations may be infinite.
For $X \in L^1$, this is equivalent to the existence of $Y' \diseq Y$ and
$X' \diseq X$ with $Y = \E (X \| \G)$ for some $\sf$ $\G$.  Also, if
$X \in L^1$ then $Y \preceq X$ is equivalent to the conjunction: 
$\E Y = \E X$ and
$\E (Y - \lambda)^+ \leq \E (X - \lambda)^+$ for all real $\lambda$.
[To see that $Y \preceq X$ implies $\E Y = \E X$, let $\phi$ be linear.
Now assuming $\E Y = \E X$, it suffices in showing $Y \preceq X$ to
consider convex $\phi$ with bounded derivative.  Since $\E Y = \E X$, 
we may add a linear term and assume $\phi$ is increasing.  Such a 
$\phi$ may be written as $\int (x - \lambda)^+ dF (\lambda)$.]
The main result of this note is as follows.  Let $\cE$ denote an
exponential random variable of mean one.

\begin{th}[discrete case] \label{discrete}
Let $Z$ be a random positive integer and $\{ \F_n \}$ be an increasing
sequence of $\sf$s.  Let 
$$Y = \sum_{n=0}^\infty \P (Z = n+1 \| \F_n) .$$
Then $Y \preceq \cE$.
\end{th}

In order both to facilitate the proof and to accommodate future applications,
I will pass to a rather general, continuous-time setting.

\begin{th}[continuous case] \label{continuous}
Let $A(t,\omega)$ be a random nondecreasing right-continuous function
with $A(0) = 0$ and $A(\infty) = 1$, and let
$\{ \F_t \}$ be an increasing, right-continuous family of $\sf$s.  
Let $\{ \AP_t \}$ be the dual previsible projection 
of $\{ A_t \}$ and let $R = \AP_\infty$.  Then $R \preceq \cE$.
\end{th}

Before proving this, let me discuss the relation between the two theorems
and the example scenarios.  It is clear how Theorem~\ref{discrete}
pertains to the insurance scenario.  To see that the upper bound
in variability is sharp, consider the continuous time insurance
problem.  Suppose that one's lifetime, $Z$, is a positive real random variable. 
The total risk is
\begin{equation} \label{eq risk}
R = \int_0^\infty \P (Z \in (t,t+dt) \| \F_t) \, ,
\end{equation}
provided the RHS makes sense.  Making sense of the RHS is where the
dual previsible projection comes in.  Let $A_t = \one_{[Z , \infty )} (t)$. 
The dual previsible projection of the increasing, right-continuous process
$\{ A_t \}$ formalizes the RHS of~(\ref{eq risk});
see~\cite[Section VI.22]{RW} for further explanation.  In the case where 
$Z$ has a density $f$ and $\F_t$ is the natural $\sf$ $\sigma (Z \wedge t)$,
this turns into the familiar
$$ R = \int_0^\infty {f(t) \over 1 - F(t)} \one_{Z > t} dt ,$$
where $F(t) = \int_0^t f(s) ds$.  It is well known that this
has a mean-one exponential distribution independent of $f$.  In fact
this is true under much more general conditions, for instance when 
$Z$ is a totally inaccessible stopping time and $\F_t$ is its natural
filtration (see~\cite[prop. 3.28]{Je}).  Two cases where the
variability is less are the extreme cases: (1) $\F_t$ is trivial
for all $t$, so $R = \int_0^\infty \P (Z \in (t,t+dt)) \equiv 1$; 
and~(2) $\F_t = \sigma (Z)$ for all $t$, in which case $R = 
\int_0^\infty d\one_{t \leq Z} \equiv 1$ again.  In general, the
insurance company will be happy to know that the variance of
the total premium of a policy based on up-to-date information 
will be less than the (easily computable) variance based on no
updated information.  

For the second scenario, let $\{ s_1 , \ldots , s_r \}$ be an ordering
of the random set $S$, where $r \leq K$ is a random variable.  For 
$1 \leq j \leq K$, let
$$Z_j = i \mbox{ if } e(i) = s_j$$
and $Z_j = \infty$ if $j > r$.  Let 
$$Y_j = \P (Z_j = \infty \| \F_\infty) + \sum_{0 \leq n < \infty} 
   \P (Z_j = n \| \F_{n-1}) .$$
It is easy to see that Theorem~\ref{discrete} extends to show that
$Y_j \preceq \cE$ for all $j$.  Let $Y_j'$ be the same as $Y_j$
but without the term $\P (Z_j = \infty \| \F_\infty)$, and let 
$W' = W \one_{|S| \leq K}$.  Then
\begin{equation} \label{eq n1}
W' \leq \sum_{j=1}^K Y_j' \leq \sum_{j=1}^K Y_j \preceq K \cE .
\end{equation}
Thus, by an easy calculation, $\P (W' > \lambda K) \leq e^{1 - \lambda}$.
(The inequality~(\ref{eq n1}) may also be derived directly from 
Theorem~\ref{continuous}.)

The pivotal bond version of the problem comes from a paper of 
H.\ Kesten on first-passage percolation, \cite{Ke}.  Here, the 
method of bounded differences (an Azuma type inequality, 
c.f.\ Wehr and Aizenman \cite{WA}) is used to bound the
variability of a first-passage time in terms of a conditional square
function that turns out to be of the form discussed above.  
Kesten~\cite[Theorem 3]{Ke} isolates the part of the argument
that requires an upper tail bound on the conditional square function.
Steps~2 and~3 of the Kesten's proof~\cite[Section 5]{Ke} may be 
replaced by the result
$$\P \left [ \sum_{k=1}^N \E (U_k \| \F_{k-1}) \geq R , \sum_{k=1}^N
   U_k \leq T \right ] \leq e^{1 - R/T} ,$$
gotten by applying Theorem~\ref{continuous} to $T^{-1} \sum U_k 
\one_{\sum U_k \leq T}$.  

Finally, to see that Theorem~1 is a special case of Theorem~2, begin with
the hypotheses of Theorem~1 and let $A_t = \one_{Z \leq t+1/2}$.
Let $\F_t' = \F_{\lfloor t \rfloor}$.  Then $\{ \F_t' \}$ 
is right-continuous and applying Theorem~2 gives
$$Y = \sum_{t + {1 \over 2} \in \Z^+} \P (Z = t+{1 \over 2} \| 
   \F_{t - {1 \over 2}}) = \AP_\infty \preceq \cE . $$

\section{Proofs}

Let $\{ \AP_t \}$ be the dual previsible projection of $\{ A_t \}$ 
as before, and let $\{\AO_t\}$ be the optional projection of $\{ A_t \}$;
the optional projection is a c\'adl\'ag process such that for each $t$,
$\AO_t$ is a version of $\E (A_t \| \F_t )$.  
\begin{lem} \label {lem 1}
The optional process $\{ M_t \}$ defined by
$$ M_t = e^{\AP_t} (1 - \AO_t) $$
is a supermartingale with respect to $\{ \F_t \}$.
\end{lem}

The intuition behind this is pretty clear: $\AP_{t+dt} - \AP_t$ is the
expected value of $\AO_{t+dt} - \AO_t$, so the total expected increase 
is $M_t (d \AP_t) - e^{\AP_t} \E (d \AO_t)$ which is never greater 
than zero.  The proof is based on the following formula:
\begin{eqnarray}
M_t - M_s & = & \int_s^t e^{\AP_{r-}} (1 - \AO_{r-}) \, d\AP_r 
   - \int_s^t e^{\AP_{r-}} \, d\AO_r \nonumber \\
&& + \sum_{s < r \leq t} \left [ e^{\AP_r} (1 - \AO_r) - e^{\AP_{r-}} 
   (1 - \AO_{r-}) \right. \nonumber \\ 
&& \left. - e^{\AP_{r-}} (1 - \AO_{r-}) (\AP_r - \AP_{r-}) +
   e^{\AP_{r-}} (\AO_r - \AO_{r-}) \right ] \; . \label{eq nn1}
\end{eqnarray}

This formula may be derived from~\cite[page 334-335]{DM} by 
the following observation: since $A$ is increasing, $\AO$
is a submartingale and $\AP$ is increasing; then by~\cite[VIII~(19.3)]{DM},
the square bracket terms in~\cite[VIII (27.1)]{DM} vanish, resulting
in~(\ref{eq nn1}).  A more direct derivation without using the 
full strength of the stochastic It\^o formula~\cite[VIII (27.1)]{DM}
is possible.  Observe for later use that $\AP_t - \AO_t$ is a
martingale: assuming without loss of generality that $A_0 = 0$,
one has for any stopping time $T$,
$$\E (\AO_T - \AP_T) = \E ( A_T - \int \one_{[0,T]} \, d\AP_s) \, ;$$
since $\one_{[0,T]}$ is left continuous and adapted, it is predictable,
and hence this becomes $\E (A_T - \int \one_{[0,T]} \, dA_s) = 0$. 

\noindent{\sc Proof of Lemma}~\ref{lem 1}: Rewrite the two integral 
terms in~(\ref{eq nn1}) as
$$ \int_s^t (- \AO_{r-}) e^{\AP_{r-}} \, d\AP_r + \int_s^t
    e^{\AP_{r-}} \, d(\AP_r - \AO_r) .$$
Using $\DAO_r$ (respectively $\DAP_r$) to denote $\AO_r - \AO_{r-}$
(respectively $\AP_r - \AP_{r-}$), rewrite the summation as 
$$ \sum_{s < r \leq t} e^{\AP_{r-}} \left [ e^{\DAP_r} (1 - \AO_r)
   - (1 - \AO_{r-}) - (1 - \AO_{r-}) \DAP_r + \DAO_r \right ] .$$
Combining the first, second and fourth terms inside the
square brackets yields 
$$(e^{\DAP_r} - 1) (1 - \AO_r) ,$$
while expanding the third term yields 
$$ - (1 - \AO_r) \DAP_r - (\DAP_r)^2 + \DAP_r (\DAP_r - \DAO_r) .$$
The quantity in square brackets may therefore be rewritten as
$$(e^{\DAP_r} - 1 - \DAP_r) (1 - \AO_r) - (\DAP_r)^2 + \DAP_r
   (\DAP_r - \DAO_r) $$
and equation~(\ref{eq nn1}) now becomes
\begin{eqnarray*}
M_t - M_s & = & \int_s^t (- \AO_{r-}) e^{\AP_{r-}} \, d\AP_r + \int_s^t
    e^{\AP_{r-}} \, d(\AP_r - \AO_r) \\
&& + \sum_{s < r \leq t} e^{\AP_{r-}} (-\AO_r) (e^{\DAP_r} - 1 - \DAP_r) +
   \sum_{s < r \leq t} e^{\AP_{r-}} (e^{\DAP_r} - 1 - \DAP_r - (\DAP_r)^2) \\
&& + \sum_{s < r \leq t} e^{\AP_{r-}} \DAP_r ( \DAP_r - \DAO_r) .
\end{eqnarray*}
The conditional expectations given $\F_s$ may be seen to be nonpositive
term by term.  The first integral is everywhere nonpositive.  The
second is the integral of a previsible process against a martingale
and hence has zero expectation given $\F_s$.
The first summation is everywhere nonpositive, as is the second, 
since $\DAP_r = \E (\Delta A_r \| \F_{r-}) \in [0,1]$ for all $r$,
and $e^z \leq 1 + z + z^2$ for $z \in [0,1]$.  Finally, the third 
summation is the integral of the previsible process $e^{\AP_{r-}} \DAP_r$
against the martingale $\AP_r - \AO_r$, and therefore has
zero conditional expectation given $\F_s$.  Thus $M_t$ is 
a supermartingale.

\noindent{\sc Proof of Theorem}~\ref{continuous}:
Fix a real $\lambda > 0$ and define a stopping time $\tau = 
\inf \{ t \geq 0 : \AP_t \geq \lambda \}$.  The purpose of the argument
between here and~(\ref{eq n2}) is to handle the case where,
due to jumps, $\AP_{\tau -} < \lambda$.  If you are not worried
about jumps, skip ahead to~(\ref{eq n2}) and read only the first
expression inside each subsequent expectation.

Since $\tau$ is previsible, there are times $\tau_n \neq \tau$ 
increasing to $\tau$ almost surely, and it follows that 
$$\E M_{\tau -} = \E \lim M_{\tau_n} \leq \liminf \E M_{\tau_n} \leq 
   \E M_0 .$$
Now define a random variable $X$ by setting $X=0$ when $\tau = \infty$,
setting $X = e^{\lambda} (1 - \AO_\tau)$ when $\DAP_\tau = 0$,
and otherwise setting
$$ X = e^\lambda \left [ {\AP_\tau - \lambda \over \AP_\tau - \AP_{\tau-}}
   (1 - \AO_{\tau-}) +  {\lambda - \AP_{\tau-} \over \AP_\tau - \AP_{\tau-}}
   (1 - \AO_\tau) \right ] .$$
The following computation shows that $\E X < \E M_{\tau-} \leq 1$ in
the case where $\DAP_\tau \neq 0$.  Make use of the facts that 
$\AP_\tau, \AP_{\tau-}$ and $\AO_{\tau-}$ are all in $\F_{\tau-}$,
and $\E (\DAO_\tau - \DAP_\tau \| \F_{\tau-}) = 0$ to write:
\begin{eqnarray*}
&& \E (X - M_{\tau-} \| \F_{\tau - }) \\[2ex]
& = & (e^\lambda - e^{\AP_{\tau - }}) (1 - \AO_{\tau-}) - e^\lambda
   (\lambda - \AP_{\tau - }) \\[2ex]
& \leq & e^\lambda - e^{\AP_{\tau -}} - e^\lambda (\lambda - \AP_{\tau - }) .
\end{eqnarray*}
This is less than or equal to 0 since $e^z - e^y - e^z (z-y) \leq 0$
for $z \geq y \geq 0$.  In the cases $\DAP_\tau = 0$ or $\tau = \infty$,
the conclusion that $\E (X - M_{\tau -} \| \F_{\tau -} ) \leq 0$ still
holds, some terms having dropped out of the above computation.

Combining this result with the fact that $\{ M_t \}$ is a supermartingale
shows that $e^{-\lambda} = e^{-\lambda} \E M_0 \geq e^{-\lambda} 
\E M_\tau \geq e^{-\lambda} \E X$.  Thus
\begin{equation} \label{eq n2}
e^{-\lambda} \geq \E \left [ (1 - \AO_\tau) +  (\AO_\tau - \AO_{\tau-}) 
   {\AP_\tau - \lambda \over \AP_\tau - \AP_{\tau - }} \right ] .
\end{equation}
Taking conditional expectations with respect to $\F_{\tau-}$ shows
that $\AO$ may be replaced by $A$, yielding
$$ e^{-\lambda} \geq \E \left [ (1 - A_\tau) +  (A_\tau - A_{\tau-}) 
   {\AP_\tau - \lambda \over \AP_\tau - \AP_{\tau - }} \right ] .$$
Since $A_\infty = 1$, we may write the RHS as the stochastic integral
$$\E \int \left ( \one_{(\tau , \infty)} (t) + {\AP_\tau - \lambda \over 
   \AP_\tau - \AP_{\tau - }} \one_{t = \tau} \right ) \, d A_t .$$
The integrand is previsible, so this becomes
\begin{eqnarray*}
e^{-\lambda} & \geq & \E \int \left ( \one_{(\tau , \infty)} (t) + {\AP_\tau - 
   \lambda \over \AP_\tau - \AP_{\tau - }} \one_{t = \tau} \right ) \, 
   d \AP_t \\[2ex]
& = & \E (\AP_\infty - \lambda)^+ .
\end{eqnarray*}
Thus the total risk $R$ satisfies $\E (R - \lambda)^+ \leq e^{-\lambda}$
for every positive $\lambda$, which, along with the fact that $\E R = 1$,
suffices to prove $R \preceq \cE$.   $\Cox$

\end{document}